\documentclass[11pt]{amsart}
\usepackage{amsmath,amssymb}
\usepackage{graphicx,color}
\usepackage[all]{xy}

\newtheorem{thm}{Theorem}[section]
 
\newtheorem{lem}[thm]{Lemma} 

\newtheorem{prop}[thm]{Proposition} 
\newtheorem{rem}[thm]{Remark} 
\newtheorem{exam}[thm]{Example}

\newcommand{\A}{\mathcal{A}}

\newcommand{\Ost}{\mathcal{O}}
\newcommand{\Proj}{\mathbb{P}}
\newcommand{\C}{\mathbb{C}}

\newcommand{\Q}{\mathbb{Q}}

\newcommand\codim{\textup{\textrm{codim}\,}}

\def\SM{{\rm SM}}

\def\bz{\mbox{\boldmath $z$}}

\def\b0{\mbox{\boldmath $0$}}

\def\imageconst{\alpha_{\mbox{\rm \tiny im}}}

\begin{document}
\title[Projective surfaces and $3$-folds with ordinary singularities]
{Classical formulae on 
projective surfaces and $3$-folds with ordinary singularities, revisited}
\author[T.~Sasajima]{Takahisa Sasajima}
\address[T.~Sasajima]{Nakagyou-ku Koromono-tana-dori Oike-agaru Shimo-myoukaku-ji 192-3, 
Kyoto 604-0024, Japan}
\email{sasajidegozaimath@gmail.com}
\author[T.~Ohmoto]{Toru Ohmoto}
\address[T.~Ohmoto]{Department of Mathematics, 
Graduate School of Science,  Hokkaido University,
Sapporo 060-0810, Japan}
\email{ohmoto@math.sci.hokudai.ac.jp}
\dedicatory{
Dedicated to Professor Takuo Fukuda 
 on the occasion of his 75th birthday.}
%
%
%
\maketitle

\begin{abstract} 
From a viewpoint of global singularity theory, 
we revisit classical formulae of Salmon-Cayley-Zeuthen 
on numerical characters of projective surfaces and 
analogous formulae of Severi-Segre-Roth on those of $3$-folds. 
The theory of universal polynomials associated to local and multi-singularities of maps 
provides a unified way for reproducing those classical enumerative formulae. 
\end{abstract}

 
\section{Introduction} 
In this paper, from a viewpoint of global singularity theory, 
we revisit classical formulae of Salmon-Cayley-Zeuthen for 
projective surfaces in $3$-space 
and analogous formulae of Severi-Segre-Roth for $3$-folds in $4$-space.  
For simplicity, 
we work over $\C$ throughout. 
In classical literature, 
numerical projective characters of a subvariety of projective space 
are defined 
by the degree of individual singularity loci,  
e.g. the number of triple points, 
and also by the degree of {\it polars}, i.e. 
the degree of the critical loci of linear projections. 
For a surface with {\it ordinary singularities} in $\Proj^3$ 
(see \S 3.1 for the definition), there are known five relations 
among nine numerical characters, which were originally given by 
Salmon \cite[XVII]{Salmon} 
(cf. Baker \cite[IV]{Baker} 
and Semple-Roth \cite[IX, \S 3]{SR}); 
those have rigorously been reformulated within modern algebraic geometry  
by Piene \cite{Piene}. 
For a $3$-fold in $\Proj^4$, called {\it primals} classically, 
there are analogous works of Roth  \cite{Roth} (and also  Severi, B. Segre etc);  
recently, new formulae on Chern numbers 
have been shown by Tsuboi \cite{Tsuboi05, Tsuboi09} 
using modern intersection theory. 

In contrast, 
our approach is based on 
singularity theory of holomorphic maps and complex cobordism theory. 
The aim of this paper is to demonstrate 
an effective way to 
reproduce and generalize those classical formulae 
by means of universal polynomials 
associated to local and multi-singularities of maps, called {\it Thom polynomials}  
(cf. Kleiman \cite{Kleiman76, Kleiman81}, Colley \cite{Colley}, 
Feh\'er-Rim\'anyi \cite{FR04, FR12}, Kazarian \cite{Kaz03, Kaz08}, 
Ohmoto \cite{Ohmoto06, Ohmoto16}). 
For a surface with ordinary singularities in $\Proj^3$, 
there are essentially four numerical characters, which actually correspond to Chern monomials 
$1, c_1, c_1^2, c_2$ of 
the normalization of the surface.  
As an analogy, 
for a $3$-fold in $\Proj^4$ 
with ordinary singularities (see \S 4.1), 
Roth  \cite{Roth} claimed that there are seven independent characters 
among more than 20 numerical characters. 
We show that all numerical characters can explicitly be computed  
by using universal polynomials, and especially,  
independent characters correspond to seven Chern monomials 
$c_1^i c_2^j c_3^k\; (0\le i+j+k \le 3)$ of its normalization (Theorem \ref{thm}). 
In particular, 
our computation recovers exactly the same formulae as those of Tsuboi 
\cite{Tsuboi09}, although methods are completely different (Remark \ref{tsuboi}). 
It should be noted that we make use of 
an advanced version of universal polynomial  
for computing {\it weighted Euler characteristics} of singularity loci of prescribed type,  
which is based on 
equivariant Chern-Schwartz-MacPherson (CSM) classes established by the second author 
\cite{Ohmoto06, Ohmoto16}. 
Indeed, the CSM class satisfies a variant of Grothendieck-Riemann-Roch formula, so 
the additivity and the covariant functoriality of CSM classes 
are fit for interpreting or generalizing standard classical arguments 
for counting invariants, e.g. addition-deletion principle. 

A generalization to higher dimensional case would be possible. 
In fact, for a projective variety $X=f(M) \subset \Proj^n$ given as 
the image variety of a locally stable map $f: M^m \to \Proj^n$ ($m<n$), 
the composition of $f$ with a generic linear projection onto lower dimensional projective space 
becomes a locally stable map, that has been proven by 
Mather \cite{Mather73} and  Bruce-Kirk \cite{BruceKirk}. 
Therefore,  in general, numerical characters of such a subvariety $X$ would be expressed in terms of 
degrees of Chern monomials of $M$ and the hyperplane class by using 
universal polynomials. 
For another kind of classical enumerative problems, 
e.g.  counting lines which have a prescribed contact with a given surface 
(cf. Salmon \cite{Salmon}), 
universal polynomials for {\it unstable} singularities are effective, 
that has been discussed in our another paper \cite{SO}.

\section{Singularities of maps and Thom polynomials} 

\subsection{Local and multi-singularities of maps} 
All maps are assumed to be holomorphic throughout. 
Let $\kappa$ denote the relative codimension $n-m$ 
for a map $f: M \to N$ from an $m$-fold to an $n$-fold. 
As local classification of maps, 
it is very natural to think of the equivalence relation of map-germs 
via local coordinate changes of the source and the target: 
$f, g: (\C^m, 0) \to (\C^n, 0)$ are {\it $\A$-equivalent} 
if there are biholomorphic germs $\sigma:(\C^m,0) \to (\C^m,0)$ 
and $\tau: (\C^n,0) \to (\C^n,0)$ so that $g=\tau \circ f \circ \sigma^{-1}$.  
Obviously,  the equivalence is defined for any germs $f: (M, p) \to (N, f(p))$, 
and an equivalence class is called a local singularity type or $\A$-singularity type.  
We say that a map-germ $f$ is {\it stable} if 
any deformation of $f$ is trivial up to  (parametric) $\A$-equivalence. 
Stable-germs $f:(\C^m,0) \to (\C^n,0)$ are characterized by its local algebra 
$$Q_f:=\Ost_{\C^m, 0}/f^*\frak{m}_n\Ost_{\C^m, 0},$$ 
where $\Ost_{\C^m, 0}$ is the ring of holomorphic function-germs at $0$ 
and $\frak{m}_n$ is the maximal ideal of $\Ost_{\C^n, 0}$ -- 
a theorem of J. Mather  (cf. \cite{Mather78}) says 
that if stable-germs $f, g:(\C^m,0) \to (\C^n,0)$ 
have isomorphic local algebras, then they are $\A$-equivalent. 
For stable-germs, $\kappa=n-m$ is more essential rather than $m, n$. 

\begin{exam}\upshape  
Let $\kappa \ge 0$. The $A_\mu$-singularity is characterized by 
$$Q_{A_\mu}=\C[x]/\langle x^{\mu+1}\rangle.$$
Normal forms of stable $A_\mu$-singularity types 
in low dimensions are as follows: 
$$
\begin{array}{ll}
\C^2 \to \C^2& A_1: (x^2, y), \; \;\; A_2: (x^3+yx, y); \\
\C^2 \to \C^3& A_1: (x^2, xy, y); \\
\C^3 \to \C^3& A_1: (x^2, y,z), \; A_2: (x^3+yx, y,z), \; A_3: (x^4+yx^2+zx, y,z); \\
\C^3 \to \C^4& A_1: (x^2, xy, y, z). 
\end{array}
$$
Those germs are classically known \cite{Salmon}; 
the image (or the critical value locus) is a typical singularity of 
singular surfaces and $3$-folds, which will be discussed later (\S 3 and \S 4). 
\end{exam}

A {\it  multi-singularity type} is an $\A$-equivalence class of 
germs (multi-germs) $f: (M, S) \to (N, q)$ 
of holomorphic maps at finite subsets 
$S=\{p_1, \cdots, p_k\} \subset f^{-1}(q)$ in $M$. 
A multi-singularity of map-germs is  {\it stable} 
if it is stable under any deformation up to $\A$-equivalence. 
In particular, if $f: (M, S) \to (N, q)$ is a stable multi-singularty, then 
each mono-germ $f: (M, p_i) \to (N, q)$ is stable and 
the sum of their codimensions does not exceed the target dimension $n$.  
A multi-singularity type may be regarded just as a collection of 
$\A$-types $\eta_i$ of mono-germs $(\C^m, 0) \to (\C^n, 0)$; 
an {\it ordered} type is denoted 
by $\underline{\eta}:=(\eta_1, \eta_2, \cdots, \eta_r)$, 
which is used below to distinguish the first entry $\eta_1$ from others. 

A proper holomorphic map $f: M \to N$ between complex manifolds 
is {\it locally stable} if for any $y \in N$ and for any finite subset $S$ of $f^{-1}(y)$, 
the germ $f: (M, S) \to (N, y)$ is stable. 
Let $M$ and $N$ be of dimension $m$ and $n$, respectively, 
and $f$ a locally stable map. 
Given an $\A$-type $\eta$ of $(\C^m, 0) \to (\C^n, 0)$, 
we set 
$$\eta(f) = \{\; p \in M \; | \; \mbox{the germ of $f$ at $p$ is of type $\eta$}\; \}$$
and call its closure the {\it $\eta$-singularity locus} of $f$, 
which becomes an analytic closed subset of $M$. 
For a stable multi-singularity type $\underline{\eta}$,  we set 
$$\underline{\eta}(f)
:=\left\{ \; p_1 \in \eta_1(f) \; \Bigl| \;
\begin{array}{l} 
\mbox{$\exists$ distinct $p_2, \cdots, p_r \in f^{-1}f(p_1)-\{p_1\}$} \\
\mbox{s.t.  $f$ at $p_i$ is of type $\eta_i$ } 
\end{array}  \right\}$$ 
and call its closure $\overline{\underline{\eta}(f)} \subset M$  
the {\it multi-singularity locus of type $\underline{\eta}$ in the source}. 
In case of $m \le n$, the restriction map 
$$f: \overline{\underline{\eta}(f)} \to \overline{f(\underline{\eta}(f)})$$ 
 is finite-to-one; 
 let $\deg_1 \underline{\eta}$ denote 
the degree of this map, that is equal to  
the number of $\eta_1$-type appearing in the tuple $\underline{\eta}$. 
e.g. $\deg_1A_0^3=3$. 

\begin{exam} \upshape 
In case of  $(m,n)=(3,4)$,  
local stable singularity type is only $A_1$, which yields  
the critical locus $\overline{A_1(f)}$ of a stable map 
$f: M^3 \to N^4$ (it becomes a smooth curve). 
Stable multi-singularity types are the following types whose components are mutually transverse: 
\begin{flushleft}
$\;\;$ $A_0^2:=A_0A_0$ (double point locus in the source space); \\
$\;\;$ $A_0^3:=A_0A_0A_0$ (triple point curve); \\
$\;\;$ $A_0^4:=A_0A_0A_0A_0$ (quadruple points);  \\
$\;\;$ $A_0A_1$ and $A_1A_0$ (preimage of stationary points). 
\end{flushleft}
\end{exam}

\subsection{Thom polynomials}\label{section_tp}

To each local singularity type $\eta$ of map-germs $\C^m,0 \to \C^{m+\kappa},0$, 
one can assign a unique universal polynomial $tp(\eta)$ of quotient Chern classes 
$c_i=c_i(f)$ which expresses the $\eta$-singularity locus of 
appropriately generic maps $f: M \to N$  
(cf. \cite{FR04, FR12, Kaz03, Ohmoto16}): 
$$ [\overline{\eta(f)}] =tp(\eta)(f)   \in H^*(M).$$
Here $c_i(f)$ is defined by the $i$-th component of $c(f)=1+c_1(f)+c_2(f)+\cdots$ 
where  
$$c(f)=c(f^*TN-TM)=\frac{1+f^*c_1(TN)+\cdots}{1+c_1(TM)+\cdots}.$$
This fact originates in Ren\'e Thom (perhaps, his Strasbourg Lecture in 1957), 
thus  $tp(\eta)$ is usually called the \textit{Thom polynomial} of $\eta$. 
As simplest examples,  $tp(A_1)$, $tp(A_2)$ and $tp(A_3)$ in cases of $\kappa=0, 1$ 
are shown in Tables \ref{tp_codim0} and \ref{tp_codim1} below.

Furthermore, one may expect the existence of universal polynomials 
for stable multi-singularity types, which generalize 
the {\it multiple point formulae} studied by Kleiman \cite{Kleiman76, Kleiman81}. 
However the general problem is known to be very hard and 
a full conjecture has been proposed in Kazarian \cite{Kaz03}. 
Nevertheless, there are some computational results in a special case of 
maps having only $A_\mu$-singularities with $\kappa \ge 0$, 
called {\it curvilinear maps}  \cite{Kleiman76, Kleiman81, Colley, Kaz03, Kaz08}.  
An algorithm given by Colley \cite{Colley} verifies the following claim 
in at least low codimensional case (cf. \cite{Kaz08});  
To each multi-singularity type $A_{\bullet}=(A_{\mu_1}, A_{\mu_2}, \cdots)$ 
(in other words, a type of {\it multiple stationary points}) of maps with codimension $\kappa$, 
one can associate a universal polynomial $tp(A_{\bullet})$ of variables 
$c_i=c_i(f)$ and $s_I=f^*s_I(f)$ with rational coefficients such that 
it holds that 
$$
\left[\overline{A_{\bullet}(f)}\right] = tp(A_{\bullet})(f) \in H^*(M; \Q)
$$
for  proper stable curvilinear maps $f: M^m \to N^{m+\kappa}$. 
Here $s_I(f)$ denote   
the {\it Landweber-Novikov classes} of $f$ multi-indexed 
by $I=(i_1 i_2\cdots)$: 
$$s_I=s_I(f)=f_*(c_1(f)^{i_1}c_2(f)^{i_2}\cdots) \; \in H^*(N).$$
In case of $\kappa=0$ and $1$, 
the universal polynomials $tp(A_{\bullet})$ with low codimension are explicitly given  
as in Tables \ref{tp_codim0} and \ref{tp_codim1} \cite{Colley, Kaz08}. 
We often write $s_I(f)$ for 
its pullback $f^*s_I(f) \in H^*(M)$, i.e. drop the notation $f^*$ 
for short.

\begin{table}
$$
\begin{array}{l  | c | l}
\hline 
\mbox{type} & \codim & tp\\
\hline 
A_1 & 1 & c_1 \\
A_2 & 2 & c_1^2 + c_2 \\
A_1^2 & 2 & c_1s_1-4c_1^2-2c_2\\
A_3 & 3 & c_1^3+3c_1c_2+2c_3\\
A_1^3 & 3 & \frac{1}{2}\left(
\begin{array}{l}
c_1s_1^2 - 4c_2s_1 - 4c_1s_2 - 2c_1s_{01} - 8c_1^2s_1\\
 +40c_1^3+56c_1c_2+24c_3
 \end{array}
 \right)\\
A_1A_2 & 3 &  c_1s_2+c_1s_{01}-6c_1^3-12c_1c_2-6c_3\\
A_2A_1 & 3 & c_1^2s_1+c_2s_1 -6c_1^3-12c_1c_2-6c_3\\
\hline
\end{array}
$$
\caption{\small $\kappa=0$ }
\label{tp_codim0}
\end{table}
\begin{table}
$$
\begin{array}{l  | c | l}
\hline 
\mbox{type} & \codim & tp\\
\hline 
A_0^2 & 1 & s_0-c_1 \\
A_1 & 2 & c_2 \\
A_0^3 & 2 & \frac{1}{2}(s_{0}^2-s_1-2s_0c_1+2c_1^2+2c_2)\\
A_0A_1 & 3 & s_{01}-2c_1c_2-2c_3\\
A_1A_0 & 3 &  s_0c_2-2c_1c_2-2c_3 \\
A_0^4 & 3 &  
\frac{1}{3!}\left(
\begin{array}{l}
s_0^3-3s_0s_1+2s_2+2s_{01}-3s_0^2c_1+3s_1c_1\\
+6s_0c_1^2+6s_0c_2-6c_1^3-18c_1c_2-12c_3
\end{array}
\right)\\
\hline
\end{array}
$$
\caption{\small $\kappa=1$ }
\label{tp_codim1}
\end{table}

\subsection{Euler characteristics}
For a complex analytic variety $X$ which admits singular points, 
i.e. $\Ost_{X, p} \not\simeq \Ost_{\C^m, 0}$, 
the Chern class $c(TX)$ no longer exists, 
because of the lack of the tangent bundle. 
However, there is a well-behaved alternative  in homology, 
$c^{\SM}(X) \in H_*(X)$, called 
the {\it Chern-Schwartz-MacPherson class} (CSM class) of $X$  \cite{Mac}: 
the CSM class satisfies a nice functorial property and the normalization condition 
that $c^{\SM}(X)=c(TX)\frown [X]$ if $X$ is non-singular. 
In particular, 
if $X$ is irreducible, proper and possibly singular,  
the degree of the class is the (topological) Euler characteristic $\chi(X)$ 
and the top component is the fundamental class, i.e. 
$$c^{\SM}(X)=\chi(X)[pt]+\cdots +[X] \in H_*(X).$$

For a closed subvariety $X$ in an ambient manifold $M$,  
the {\it Segre}-SM class of the embedding $\iota: X \to M$ is defined by 
$$s^{\SM}(X, M):=\iota^*c(TM)^{-1}\frown c^{\SM}(X) \in H_*(X).$$ 
If $X$ is non-singular, the SSM class is just the inverse normal Chern class 
$c(-\nu_{M,X})\; (=c(TX)/\iota^*c(TM))$. 
We consider $\iota_*c^{\SM}(X), \, \iota_*s^{\SM}(X, M) \in H^*(M)$ 
via the Poincar\'e dual, unless specifically mentioned. 

In \cite{Ohmoto06}, 
an equivariant version of CSM classes has been established. 
As a direct corollary, 
it is proven in \cite[Thm. 4.4]{Ohmoto16} 
(also \cite[Thm. 5.5]{Ohmoto06}, \cite{Ohmoto07, PP}) that 
for a local singularity type $\eta$, 
there exists a unique universal polynomial $tp^{\SM}(\overline{\eta})$ in $c_i(f)$ so that  
$$c^{\SM}(\overline{\eta(f)}) = c(TM) \cdot tp^{\SM}(\overline{\eta})(f) 
\in H^*(M)$$
for appropriately generic $f: M \to N$.  
We call $tp^{\SM}(\overline{\eta})$ the {\it universal SSM class} of $\overline{\eta}$, 
or roughly the {\it higher Thom polynomial} of $\overline{\eta}$ (cf. \cite{Ando}). 
It follows from the above form of $c^{\SM}(X)$ that 
 the leading term of the universal polynomial is just the Thom polynomial of $\eta$, 
$$tp^{\SM}(\overline{\eta}) =tp(\eta)+h.o.t.,$$
and the degree gives  $\chi(\overline{\eta(f)})$. 
Let us see the simplest examples. For 
$$
\begin{array}{ll}
A_1: (x, \bz) \mapsto (x^2, \bz) & (\kappa=0)\\ 
A_1: (x, y, \bz) \mapsto (x^2, xy, y, \bz) & (\kappa=1)
\end{array}
$$
low degree terms of the universal SSM class are computed as follows (see \cite[p.222]{Ohmoto16}): 
\begin{equation}\label{tpsm_codim0}
tp^{\SM}(\overline{A_1})=\left\{
\begin{array}{ll}
c_1-c_1^2+c_1^3-(c_1^4+c_2^2-c_1c_3)+\cdots & (\kappa=0)\\
c_2  - (c_1 c_2 + c_3) + \cdots& (\kappa=1).
\end{array}
\right.
\end{equation}
Then, for stable maps $f: M^m \to N^{m+\kappa}$ with low dimensions and $M$ compact, 
the Euler characteristics of the locus $S(f) :=\overline{A_1(f)} \subset M$ 
is computed by  
$$
\chi(S(f))=\int_M c(TM)\cdot tp^{\SM}(\overline{A_1})(f).  
$$

Next, for multi-singularity types in low dimensions,  
it has been stated in \cite[\S 6.2]{Ohmoto16} that 
there is a universal polynomial of $c_i=c_i(f)$ and $s_I=s_I(f)$ such that 
$$c^{\SM}(\overline{A_\bullet(f)}) = c(TM) \cdot tp^{\SM}(\overline{A_\bullet})(f) 
\in H^*(M; \Q),$$ 
for proper stable curvilinear maps $f: M \to N$. If $M$ is compact, we have 
$$
\chi(\overline{A_\bullet(f)}) = \int_M c(TM)\cdot tp^{\SM}(\overline{A_\bullet})(f). 
$$
In \cite[p.248]{Ohmoto16}, low degree terms of 
some universal polynomials $tp^{\SM}(A_\bullet)$ are explicitly computed,  
e.g., 
\begin{equation}\label{tpsm_codim1}
\textstyle 
tp^{\SM} (\overline{A_0^2}) 
= (s_0 -c_1) + \frac{1}{2}(2c_2 + 2c_1s_0 - s_0^2 - s_1) +\cdots. \qquad (\kappa=1)
\end{equation}
Note that the leading term $s_0 -c_1$ is just $tp (\overline{A_0^2})$ (the double point formula).  
Moreover, 
it is also shown in Theorems 6.5 and 6.10 in \cite{Ohmoto16} 
that for proper stable curvilinear maps $f: M^m \to N^{m+1}$, 
the Euler characteristics of 
the image hypersurface $X=f(M)$ 
and the (singular) double locus $D=\overline{f(A_0^2(f))}$ in $N$ 
are expressed as 
$$\chi(X)=\int_M c(TM) 
\cdot tp^{\SM}(\imageconst)(f)$$
and 
$$\chi(D)=\int_M c(TM) 
\cdot tp^{\SM}(\imageconst(2))(f), 
$$ 
where 
\begin{eqnarray*}
tp^{\SM}(\imageconst)
&=&\textstyle 1+\frac{1}{2}(c_1-s_0)
 +\frac{1}{6}(s_0^2+2s_1-2c_1s_0-c_1^2-c_2)\\
&&\textstyle +\frac{1}{24}\left(
\begin{array}{l}
\textstyle 2c_1^3-10c_1c_2+2c_1^2s_0+2c_2s_0+3c_1s_0^2\\
\textstyle  -s_0^3+14s_{01}+5c_1s_1-5s_0s_1-6s_2 
\end{array}
\right) \\
&&+ \cdots.\\
tp^{\SM}(\imageconst(2))
&=&\textstyle
 \frac{1}{2}(s_0-c_1)+\frac{1}{6}(-c_1^2+5c_2+4c_1s_0-2s_0^2-s_1)\\
 &&  \textstyle 
+ \frac{1}{24}\left(
\begin{array}{l}2c_1^3+38c_1c_2+24c_3+2c_1^2s_0-22c_2s_0-9c_1s_0^2\\
+3s_0^3-14s_{01}-7c_1s_1+7s_0s_1+2s_2
\end{array}
\right)\\
&&+\cdots.
\end{eqnarray*}
Those universal polynomials will be used in the following sections.

\section{Classical  formulae of surfaces in $3$-space}
A simple application of Thom polynomials recovers formulae of 
Salmon-Cayley-Zeuthen of surfaces in $3$-space. 
This is a prototype for a further application to $3$-folds which will be discussed 
in the next section. 

\subsection{Surfaces with ordinary singularities}
\begin{figure}
  \includegraphics[width=10cm]{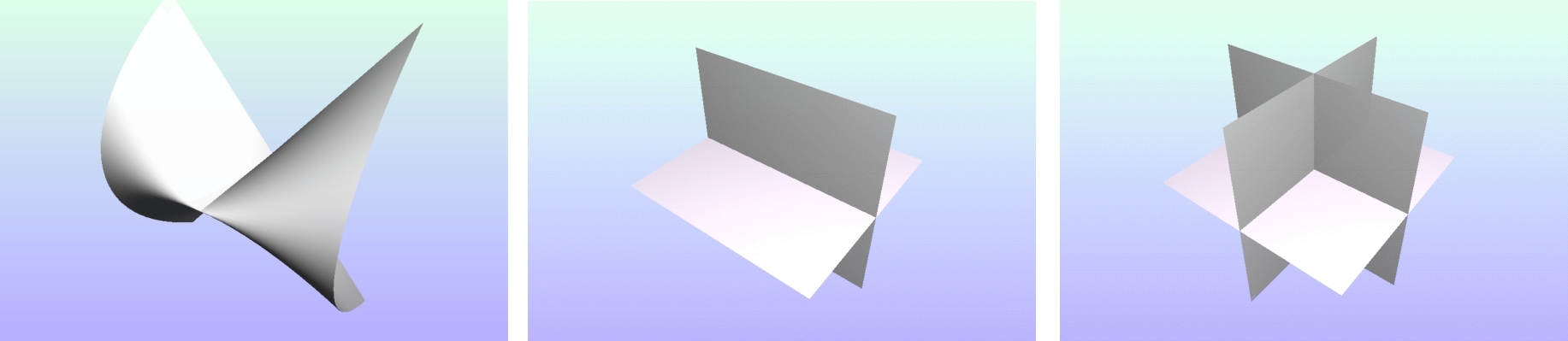}
\caption{Surface with ordinary singularities in $3$-space.}
\label{2-3}
\end{figure}

Let $X \subset \Proj^3$ be a reduced surface of degree $d$ having 
only {\it ordinary singularities}.  That is,  
the singular locus of $X$ consists of crosscap points 
and transverse double and triple points, 
which are locally defined, respectively, by equations 
up to local analytic coordinate changes (Figure \ref{2-3}): 
$$xy^2-z^2=0, \quad xy=0, \quad xyz=0.$$

Note that 
the $A_1$-singularity $\C^2 \to \C^3$, $(u, v)\mapsto (x,y,z)=(v^2, u,  uv)$, 
is a normalization of the crosscap. 
Hence, the normalization $M$ of $X$ becomes smooth, 
and the composition map of the projection and the inclusion into the ambient space, 
denoted by $f:M\rightarrow \Proj^3$ with $X=f(M)$, 
is locally stable   (cf. \cite{Piene, Mather78}).  
We do not need the condition that $M$ is embedded in some higher dimensional projective space 
and  $f$ is realized by a linear projection. 
Let $\Gamma$ be the double point locus in the source and 
$D$ its image (the double curve): 
$$\Gamma=\overline{A_0^2(f)} \subset M, 
\quad D=f(\Gamma) \subset \Proj^3.$$
Note that $\Gamma$ is an immersed curve having nodes at $A_0^3(f)$ 
(the preimage of triple points) on $M$; 
$D$ is an immersed space curve with triple points. 
The crosscap points are located on $D$, where $\Gamma \to D$ is doubly ramified. 

From now on, we always use cohomology with rational coefficients: $H^*=H^*(-, \Q)$. 
We denote the hyperplane class of $\Proj^3$ by 
$$a=c_1(\Ost_{\Proj^3}(1))\in H^2(\Proj^3).$$  
Put 
$$
f_*(1)= d  a, \;\;
f_*c_1(TM)=\xi_1  a^2, \;\;
f_*(c_1(TM)^2)=\xi_2 a^3, $$
$$
f_* c_2(TM)=\xi_{01} a^3\in H^*(\Proj^3). 
$$
Then 
\begin{align*}
&c(f)=c(f^*T\Proj^3-TM)=\frac{(1+\tilde{a})^4}{1+c_1(TM)+c_2(TM)}\\
& =1+(4\tilde{a}-c_1(TM))+(6\tilde{a}^2-4\tilde{a} c_1(TM)+c_1(TM)^2 -c_2(TM))+\cdots
\end{align*}
($\tilde{a}:=f^*a\in H^2(M)$) and 
\begin{align*}
	s_0(f)&=f_*(1)=da, \\
	s_1(f)&=f_*(c_1(f))=(4d-\xi_1)a^2,\\
	s_2(f)&=f_*(c_1(f)^2)=(16d-8\xi_1+\xi_2)a^3,\\
	s_{01}(f)&=f_*(c_2(f))=(6a-4\xi_1+\xi_2-\xi_{01})a^3.
\end{align*}

\subsection{Numerical projective characters}
Classically, 
nine numerical characters of $X$ are introduced  \cite{Baker, SR}: 
\begin{itemize}
\item 
$\mu_0=d$: the degree of $X$
\item 
$\mu_1$: the rank of $X$ (the first polar)
\item 
$\mu_2$: the class of $X$ (the second polar)
\item 
$\kappa$: the number of cusps in generic projection $X \to \Proj^2$ 
\item 
$\epsilon_0$: the degree of the double curve $D$ 
\item 
$\epsilon_1$: the rank of $D$
\item 
$\rho$: the class of immersion of $D$ in $X$
\item 
 $C$: the number of crosscaps 
\item 
$T$: the number of triple points.  
\end{itemize}

We express the above numerical characters 
in terms of $d, \xi_1, \xi_2, \xi_{01}$ 
by using Thom polynomials. 
First,  $C, T, \epsilon_0$ immediately follows from 
Thom polynomials of $A_1$, $A_0^2$ and $A_0^3$ (Table \ref{tp_codim1}) 
applied to $f: M \to \Proj^3$: 
\begin{align*}
C&\textstyle = {\displaystyle \int} f_*(Tp(A_1)(f))=6d-4\xi_1+\xi_2-\xi_{01}, \\
T&\textstyle = \frac{1}{3} {\displaystyle \int} f_*(Tp(A_0^3)(f))\\
&\textstyle =\frac{1}{6}(44d-12d^2+d^3+(3d-24)\xi_1+4\xi_2-2\xi_{01}), \\
\epsilon_0
&\textstyle 
=\frac{1}{2}{\displaystyle \int} f_*(Tp(A_0^2)(f))\cdot a=\frac{1}{2}(d^2-4d+\xi_1). 
\end{align*}
Conversely, we can express $\xi_1, \xi_2, \xi_{01}$ in terms of $\epsilon_0, C, T$: 
\begin{lem}\label{enriques} 
It holds that 
\begin{align*}
	\xi_1&=d(4-d)+2\epsilon_0,\\
	\xi_2 &= d(d-4)^2 +(16-3d)\epsilon_0 +3T -C,\\
	\xi_{01} &= d(d^2-4d+6)+(8-3d)\epsilon_0 +3T -2C.
\end{align*}
\end{lem}

\begin{rem}\upshape 
Lemma \ref{enriques} is classical  (see \cite[Prop.1]{Piene}). 
In old literature, 
the degree $\mu_0=d$, the rank $\mu_1$, the class $\mu_2$ 
and the degree of pinch points $C$ of $X$ are more frequently used; 
they are called {\it elementary numerical characters}. 
A major classical interest was to express intrinsic invariants 
of projective varieties 
in terms of their extrinsic invariants (=numerical characters). 
For instance, 
two basic invariants are introduced for projective surfaces: 
$$\omega:=\mu_2-6\mu_1+9\mu_0+C+1, 
\quad I:=\mu_2-2\mu_1+3\mu_0-4,$$ 
called the {\it Castelnuovo-Enriques invariant} 
and  the {\it Zeuthen-Segre invariant}, respectively 
(\cite[Chap.V]{Baker}, \cite[pp.223--224]{SR}). 
Their sum is an absolute birational invaraiant: 
$\omega+I=12p_a+9$, where $p_a$ is the arithmetic genus of $X$. 
That is verified in our context as follows: 
by the computation of $\mu_1$ and $\mu_2$ below, we see  
$$ \omega=\xi_2+1, \quad  I=\xi_{01}-4, 
$$
thus by using the Hirzebruch-Riemann-Roch formula, 
$$\omega+I=\int (c_1(TM)^2+c_2(TM))-3=12\chi(M, \Ost_M)-3 
=12p_a+9.$$ 
\end{rem}

To compute other five characters, 
we need to discuss generic projections of $X$ and $D$. 
Consider a generic projection of $X$ from 
a general point $p$ in $\Proj^3$. 
Put 
$$g=pr\circ f:M\rightarrow \Proj^2$$
to be the composition of $f$ with the canonical projection $pr:\Proj^3-p\rightarrow \Proj^2$. 
Obviously $g$ has the degree $d$. 
Let $S(g)$ be the closure of $A_1(g) \subset M$, i.e. 
the critical curve of the map $g$. 
In this case, the curve is non-singular and closed. 

The space $\Proj^3-p$ is isomorphic 
to the total space of the line bundle $\Ost_{\Proj^2}(1)$, 
and the isomorphism $H^*(\Proj^2)=H^*(\Proj^3-p)$ is given 
by the pullback of the projection. 
Denote by $a'$ the hyperplane class of $\Proj^2$; 
then $\tilde{a}:=g^*a'=f^*a \in H^2(M)$. 
The quotient Chern class for $g$ is 
$c(g)=(1+\tilde{a})^3c(TM)^{-1}$
and $g_*(1)=g_*([M])=d$. 
Substitute $c_i(g)$ for $c_i$ in  (higher) 
Thom polynomials,  
then one gets the degrees and 
the Euler characteristics of corresponding loci of $g$. 

We are now ready to compute the remaining five numerical characters. 

\

\noindent $\bullet\;\;$ 
\underline{The degree $\kappa$}  is the ramification divisor of $S(g)\rightarrow g(S(g))\subset \Proj^2$, 
that is just the number of cusps of $g$. Thus 
\begin{align*}
	\kappa=\int f_*(tp(A_2)(g))&=12d-9\xi_1+2\xi_2-\xi_{01}\\
	&=d(d-1)(d-2)+(6-3d)\epsilon_0 +3T.
\end{align*}

\noindent $\bullet\;\;$ 
\underline{The {\it rank} $\mu_1$ of $X$} is defined as the degree of the first polar of $X$, i.e. 
the degree of the locus of points at which tangent planes contain a given general point $p$. 
That is equal to the degree of $S(g)$: 
$$
\mu_1=\int f_*(tp(A_1)(g))=3d-\xi_1
=d(d-1)-2\epsilon_0.
$$

\noindent $\bullet\;\;$ 
\underline{The {\it class} $\mu_2$ of $X$} is defined as the degree of the second polar of $X$, i.e. 
the degree of the locus of points at which tangent planes contain a given general line $L$. 
First, consider a generic projection $g$ to $\Proj^2$ from $p$. 
The contour curve $g(S(g)) \subset \Proj^2$ has the degree $\mu_1$ and 
$\mu_2$ polar points from a general point $pt \in \Proj^2$;  
then $L$ is given by the line joining $p$ and $pt$.  
The composed map $S(g) \to  \Proj^2-pt \to \Proj^1$ 
has $\mu_2+\kappa$ critical points, thus 
by the Riemann-Hurwitz formula  (i.e. $tp(A_1) =c_1$)
$$\mu_2+\kappa=2 \mu_1 -\chi(S(g)).$$
The Euler characteristic of $S(g)$ is computed using (\ref{tpsm_codim0})  in \S 2.3:  
$$
\chi(S(g))=\int_M c(TM)\cdot tp^{\SM}(\overline{A_1})(g)
=-9d+9\xi_1-2\xi_2.
$$
Hence, 
\begin{eqnarray*}
	\mu_2 &=&3d-2\xi_1+\xi_{01}\\
	&=& d(d-1)^2+(4-3d)\epsilon_0+3T-2C.
\end{eqnarray*}

\noindent $\bullet\;\;$ 
\underline{The {\it class $\rho$ of immersion of $D$} in $X$} is defined by 
the number of points in $D$ 
at which one of the tangent plane of $X$ contains a given point $p$. 
Crosscap points should be excluded because of the lack of the tangent plane. 
In our situation, curves $\Gamma$ and $S(g)$ on $M$ are transverse to each other; 
in particular the intersection contains crosscap points. 
Hence $\rho$ is computed as follows: 
\begin{align*}
	\rho &=\Gamma\cdot S(g)-C\\
	&=-18d+3d^2+(11-d)\xi_1-2\xi_{01}+\xi_2\\
	&=(d-2)\epsilon_0-3T.
\end{align*}

\noindent $\bullet\;\;$ 
\underline{The {\it rank} $\epsilon_1$ of $D$} is defined 
by the degree of its polar, i.e., 
the number of planes tangent to $D$ 
which contain a given general line $L$. 
First we resolves  (non-transverse) triple points of $D$ 
to get smooth $D'$. 
Let $\pi : D'\rightarrow \Proj^1$ be the composed map with the projection from $L$. 
Then $\epsilon_0$ is just the degree of $\pi$, and 
$\epsilon_1$ is equal to the number of critical points of $\pi$. 
Hence, by $tp(A_1) =c_1$,  we have 
$$\epsilon_1 = \epsilon_0\chi(\Proj^1)-\chi(D').$$
Obviously,  $\chi(D')=\chi(D)+2T$. 
On the other hand, we see that 
\begin{eqnarray*}
\chi(D) 
&=& \int_M c(TM)\cdot tp^{\SM}(\imageconst(2)))(f)\\
&=& \textstyle 
\frac{1}{3}(7 d + 6 d^2 - d^3 - \frac{5}{2} \xi_{01} - 12 \xi_1 + \frac{7}{2} \xi_2)\\
&=&  \textstyle
(4-d)\epsilon_0 + T + \frac{1}{2}C.
\end{eqnarray*}
Combining all of the above, 
\begin{eqnarray*}
	\epsilon_1 
	&=&
	\textstyle 
	3d^2-21d+(13-d)\xi_1+\frac{3}{2}\xi_{01}- \frac{5}{2}\xi_2\\
	&=&
	\textstyle(d-2)\epsilon_0-3T-\frac{1}{2}C.
\end{eqnarray*}

As seen above, all numerical characters are expressed 
by $d, \xi_1, \xi_2, \xi_{01}$. 
These expressions 
immediately imply the following five basic relations among nine characters 
in \cite[Chapter IV]{Baker} (cf. \cite[Thm.2]{Piene}): 

\begin{prop}[\cite{Baker}]\label{baker}
The above nine characters of a projective surface in $\Proj^3$ with ordinary singularities 
satisfy  
$$
\begin{array}{clcl}
(i)& d(d-1)=\mu_1+2\epsilon_0, &
(ii) & \mu_1(d-2)=\kappa+\rho, \\
(iii) &  \epsilon_0(d-2)=\rho+3T,  &
(iv) & 2\rho-2\epsilon_1=C, \\
(v) &\mu_2+2C=\mu_1+\kappa. &&
\end{array}
$$
\end{prop}

\section{Classical formulae of $3$-folds in $4$-space}

\subsection{$3$-folds with ordinary singularities}
Let $X \subset \Proj^4$ be a reduced hypersurface of degree $d$ 
with {\it ordinary singularities}, i.e., singularities are locally given by 
the following equations: 
\begin{center}
\begin{tabular}{lcl}
crosscaps & ($A_1$)& $xy^2=z^2$ \\
double point locus &($A_0^2$)& $xy=0$ \\
triple point curve &($A_0^3$)& $xyz=0$ \\
quadruple points  &($A_0^4$)& $xyzw=0$ \\
stationary points &($A_0A_1$)& $(xy^2-z^2)w=0$. 
\end{tabular}
\end{center}
There are many examples obtained by generic linear projections of smooth $3$-folds 
sitting in higher dimensional projective space. 
In fact, it is shown that 
almost all projections of smooth $3$-folds into $\Proj^4$ are locally stable 
(Mather \cite{Mather73}). 

As with the case of surfaces in $\Proj^3$, 
the $A_1$-singularity $\C^3 \to \C^4$, 
$(u, v, t)\mapsto (x,y,z,w)=(v^2, u,  uv, t)$, 
locally gives the normalization of the crosscap locus (a curve on $X$). 
Therefore, a normalization $M$ of $X$ becomes to be smooth 
and we have a locally stable map $f: M \to  \Proj^4$ with $X=f(M)$ 
(cf. \cite{Tsuboi05, Tsuboi09, Mather78}). 
Note that it is not needed that $f$ is realized by some linear projection. 
The double point locus $\Gamma=\overline{A_0^2(f)}$ in the source  $M$ is 
a closed surface with ordinary singularities; 
the double curve, triple points and crosscaps of $\Gamma$ 
are just $\overline{A_0^3(f)}$, $A_0^4(f)$ and $A_0A_1(f)$, respectively. 
The critical locus $\overline{A_1(f)}$ in $M$ is an immersed curve lying on $\Gamma$, 
and denote by $C$ its image in $\Proj^4$. 
The loci $A_1A_0(f)$ and $A_0A_1(f)$ 
are just the preimage of the stationary point locus denoted by $St$. 
The image of $\Gamma$, called the {\it double surface} $ D:=f(\Gamma)$, is 
a projective surface which has singularities along 
the target triple point curve $T=f(\overline{A_0^3(f)})$. 
Note that $f|_\Gamma: \Gamma \to D$ is ramified along $C$, 
and the source triple point locus $\overline{A_0^3(f)}$ is ramified on $T$ at $St$.

\subsection{Numerical projective characters}
As a natural generalization, 
several attempts to find invariants of $3$-folds (primals) 
were intensively made 
by algebraic geometers in that age, 
e.g. B. Segre, Severi, Todd, Eger and Roth; 
there are 
(at least) twenty numerical characters of a $3$-fold $X$ in $\Proj^4$: 
 \begin{itemize}
 \item 
$d$: the degree of $X$; 
\item 
$\mu_0$: the degree of the double surface $D=f(\overline{A_0^2})\subset X \subset \Proj^4$; 
\item 
$t$: the degree of the triple point curve $T=f(\overline{A_0^3})$; 
\item 
$q$: the number of quadruple points $Q=f(A_0^4)$; 
\item 
$s_t$: the number of stationary points $St=f(A_0A_1)=f(A_1A_0)$; 
\item 
$\gamma$: the degree of the critical curve $C=f(\overline{A_1})$; 
\item 
characters associated to projections of $X, D, T$ and $C$, i.e., 
\begin{itemize}
\item  
 rank $m_1$,  first class $m_2$ and  class $m_3$ of $X$; 
 \item 
 rank $\mu_1$ and  class $\mu_2$ of  $D$; 
\item 
 rank of $T$ and $C$; 
\item 
 class of immersions of $T$ and $C$ in $X$; 
\item 
 class of immersions of $T$ and $C$ in $D$; 
\item 
apparent characters of $D$: 
crosscaps, triple points and the double curve 
for the image of a generic projection $D \to \Proj^3$. 
\end{itemize}
\end{itemize}
Roth  \cite{Roth} claimed that there are seven independent characters among them. 
However,  it is not easy to follow his arguments because 
several intermediate characters are newly introduced and 
a number of relations among many those quantities are discussed together. 
Therefore, 
instead of tracking Roth's proofs, 
we show directly that 
there are seven independent characters corresponding to 
Chern monomials of $M$, and all other characters can be expressed by those ones. 
Put 
$$f_*(1)= d  a, \;\; f_*(c_1(TM))=\xi_1  a^2, \;\; \cdots,  \;\; f_*(c_3(TM))
=\xi_{001} a^4$$
in $H^*(\Proj^4; \Q)=\Q[a]/\langle a^5\rangle$. 
The rank of $C$ is written by its degree $\gamma$ and $\chi_C:=\int_C c_1(TC)$, 
so we use $\chi_C$, instead. 

\begin{thm}\label{thm}
Any numerical projective characters of a $3$-fold in $4$-space 
with ordinary singularities 
are expressed in terms of $d, \xi_1, \cdots, \xi_{001}$. 
In particular, 
all numerical characters are generated by  
$d, \mu_0, t, q, s_t, \gamma$ and $\chi_C$. 
\end{thm}

Theorem \ref{thm} follows from 
Examples \ref{lem_chi},  \ref{classX} and \ref{lem_classD}  below. 
We conjecture that 
the assertion would be a general phenomenon in any dimension. 

\begin{exam}  \label{lem_chi} 
\upshape 
Seven quantities $d, \mu_0, t, q, s_t, \gamma, \chi_C$ 
are expressed in terms of $d, \xi_1, \cdots, \xi_{001}$, 
and vice versa, as in  Table \ref{tp_computation} and Table \ref{chern_number}. 
In fact, they are immediately computed by
Thom polynomials in Table \ref{tp_codim1} ($\kappa=1$)  
and $tp^{SM}(\overline{A_1})$ as in (\ref{tpsm_codim0}) in \S 2.3    
applied to our stable map $f: M^3 \to \Proj^4$, 
and conversely, the degree of Chern monomials of $M$ are solved. 

\begin{table}
\begin{eqnarray*}
\mu_0 &=&\textstyle 
\frac{1}{2} {\displaystyle \int} f_*(tp(A_0^2)(f))\cdot a^2
=\textstyle
\frac{1}{2}(-5 d + d^2 + \xi_1),\\
t  &=&\textstyle 
\frac{1}{3} {\displaystyle \int} f_*(tp(A_0^3)(f))\cdot a
\\
&=&\textstyle
\frac{1}{3} (35 d - \frac{15}{2}d^2 + \frac{1}{2}d^3 - \xi_{01} - 15  \xi_1 + \frac{3}{2} d \xi_1 + 2 \xi_2 ),\\
\gamma  &=& \textstyle
{\displaystyle \int} f_*(tp(A_1)(f)) \cdot a 
=\textstyle
10 d -  \xi_{01} - 5 \xi_1 +  \xi_2, \\
q &=& \textstyle
\frac{1}{4} {\displaystyle \int} f_*(tp(A_0^4)(f))\\
&=&\textstyle
\frac{1}{4}(
-295 d + \frac{355}{6} d^2 -5 d^3 + \frac{1}{6} d^4 + 2 \xi_{001} + 
(25  - \frac{4}{3} d) \xi_{01} \\
&&\textstyle
+ (200 - 25 d+d^2)  \xi_1 + \frac{1}{2} \xi_1^2 - 7 \xi_{11} +(- 55 + \frac{8}{3} d) \xi_2 + 6 \xi_3),\\
s_t  &=&  \textstyle
{\displaystyle \int} f_*(tp(A_0A_1)(f))\\
&=&\textstyle
-120  d + 10  d^2 + 2 \xi_{001} + (20-d) \xi_{01}  + (90 -5d)\xi_1  - 6 \xi_{11}\\
&&\textstyle
+(d-30) \xi_2  + 4 \xi_3, \\
\chi_C  &=& \textstyle
{\displaystyle \int_M} c(TM) \cdot tp^{\SM}(\overline{A_1})(f)\\
&=& \textstyle
-60 d +  \xi_{001} + 10\xi_{01} + 55 \xi_{1} - 4  \xi_{11} - 20  \xi_{2} + 3  \xi_{3}.  \\
\end{eqnarray*}
\caption{Thom polynomials applied to $f: M \to \Proj^4$}
\label{tp_computation}
\end{table}

\begin{table}
\begin{eqnarray*}
\xi_0 &=& \int f_*(1)\cdot a^3= d,\\
\xi_1&=&  \int f_*(c_1(TM)) \cdot a^2 = 5 d-d^2+2\mu_0,\\
\xi_2 &=&  \int f_*(c_1(TM)^2) \cdot a\\
&=&\textstyle 
25 d-10 d^2+d^3+(20-3 d)\mu_0+3 t-\gamma,\\
\xi_{01} &=&  \int f_*(c_2(TM))\cdot a \\
&=&
\textstyle
10 d - 5 d^2 + d^3 + (10 - 3 d) \mu_0 + 3 t - 2 \gamma,\\
\xi_3 &=& \int_M c_1(TM)^3 \\
&=& \textstyle  
125 d - 75 d^2 + 15 d^3 - d^4 + (150  - 45 d  + 4 d^2  - 2 \mu_0)\mu_0\\
&&\textstyle +  4 q - \frac{1}{2}s_t + (45  - 4 d) t +(- 10 + \frac{1}{2}d) \gamma - \chi_C, \\
\xi_{11} &=&  \int_M c_1(TM)c_2(TM)\\
&=&\textstyle  
50 d - 35 d^2 + 10 d^3 - d^4 + (70  - 30 d  + 4 d^2  - 2 \mu_0)\mu_0\\
\textstyle
&& +  4 q + (30  - 4 d) t - 5 \gamma - 2 \chi_C,\\
\xi_{001} &=&  \int_M c_3(TM) \\
&=& \textstyle
10 d - 10 d^2 + 5 d^3 - d^4 + (20  - 15 d  + 4 d^2  - 2 \mu_0)\mu_0 + 4 q\\
&&\textstyle
 + \frac{3}{2} s_t + (15  - 4 d) t + (10  - \frac{3}{2} d)\gamma  - 4 \chi_C. \\
\end{eqnarray*}
\caption{Chern numbers of $M$ and other degrees.}\label{chern_number}
\end{table}

\end{exam}

\begin{exam} \upshape \label{classX} 
We compute explicitly elementary characters of $X$ 
(degree $d$, rank $m_1$, first class $m_2$ and class $m_3$)  in a similar way as in \S 3.2. 
Let $m_0=d$, the degree of $X$. 

\

\noindent $\bullet\;\;$ 
\underline{The {\it rank} $m_1$ of $X$} is defined by  
the degree of the locus $S_1$ consisting of points 
at which tangent planes contain a given general point $p$. 
Let $g: M \to \Proj^3$ be
the composition of the normalization map $f$ with generic projection from $p$. 
Then $g$ is stable and $S_1=\overline{A_1(g)}$ is a smooth surface;  
denote by $i: S_1 \hookrightarrow M$ the inclusion. 
For $i_*(1)=tp(A_1)(g)=c_1(g)$,  we have 
$$
m_1=\int f_*(tp(A_1)(g))=4d-\xi_1= d(d-1)-2\mu_0.
$$
This is called the Cayley formula in \cite{Roth}. 
Furthermore, using $tp^{\SM}(A_1)(g)=c_1(g)-c_1(g)^2+\cdots$ in (\ref{tpsm_codim0}) 
or using an embedded resolution of $S_1$ (for computing $c_1(TS_1)^2$), 
we can see that 
\begin{eqnarray*}
&&i_*(c_1(TS_1))=-c_1(g)^2+c_1(g)c_1(TM),\\
&&i_*(c_2(TS_1))=c_1(g)^3-c_1(g)^2c_1(TM)+c_1(g)c_2(TM), \\
&&i_*(c_1(TS_1)^2)=c_1(g)^3-2c_1(g)^2c_1(TM)+c_1(g)c_1(TM)^2. 
\end{eqnarray*}

\noindent $\bullet\;\;$ 
\underline{The {\it first class} $m_2$ of $X$} is defined by  
the degree of the locus $S_2$ consisting of points 
at which tangent planes contain a given general line $L$. 
Assume $p \in L$.  
Project the image singular surface $g(S_1) \subset \Proj^3$ 
from $pt =L \cap \Proj^3$, and 
denote the composed map by $h: S_1 \to \Proj^2$; 
we can assume that $h$ is stable. 
Then $S(h):=\overline{A_1(h)}$ is a smooth curve in $S_1$; 
denote by $j: S(h) \hookrightarrow S_1$ the inclusion. 
By using $tp^{\SM}(A_1)(h)$, it is easy to see that 
$$j_*(1)=c_1(h), \quad 
j_*(c_1(TS(h)))=c_1(h)c_1(TS_1)-c_1(h)^2 \;\; \in H^*(S_1).$$  
Notice that 
the critical locus $S(h)$ contains the cuspidal points of $g$; 
in fact, 
$$S(h)=S_2 \sqcup \overline{A_2(g)} \;\; \mbox{ (disjoint)}. $$
Hence, the degree $m_2$ of $S_2$ is given by 
using $tp(A_1)=c_1$ and $tp(A_2)=c_1^2+c_2$: 
\begin{align*}
	m_2&=\int f_*i_*(tp(A_1)(h)) - \int f_*(tp(A_2)(g))\\
	&=6d-3\xi_1-\xi_{01}\\
	&=d(d-1)^2+(4-3d)\mu_0+3t-2\gamma. 
\end{align*}

\noindent $\bullet\;\;$ 
\underline{The {\it class} $m_3$ of $X$} is defined by  
the number of points at which tangent planes contain a given general plane $\Pi$. 
Assume $L \subset \Pi$. 
Now the image plane curve $h(S(h)) \subset \Proj^2$ 
admits only cusps and nodes as its singularities. 
Project it from $pt_2=\Pi \cap \Proj^2$ to $\Proj^1$;  
denote the resulting map  by $h': S(h) \to \Proj^1$. 
Also denote the projection of the cuspidal locus of $g$ by 
$$h'':=h'|_{\overline{A_2(g)}}: \overline{A_2(g)} \to \Proj^1.$$ 
Then the number of critical points of $h'$ is the sum of four quantities; 
\begin{itemize}
\item[-] 
$A=m_3=$ $\#$ polar points of $h(S_2)$ from $pt_2$, 
\item[-] 
$B=$ $\#$ cuspidal points of $h(S_2)$, 
\item[-] 
$C=$ $\#$ polar points of $h(A_2(g))$ from $pt_2$, 
\item[-] 
$D=$ $\#$ cuspidal points of $h(A_2(g))$ $=\#$ swallowtail $A_3(g)$ of $g$. 
\end{itemize}
Also the number of critical points of $h''$ is just $C+D$, and 
the number of cusps of $h$ is $B+D$. 
By using Thom polynomials
\begin{itemize}
\item[-] 
$A+B+C+D=\int tp(A_1)(h')=\int c_1(h')$; 
\item[-] 
$C+D=\int tp(A_1)(h'')=\int c_1(h'')$; 
\item[-] 
$B+D=\int tp(A_2)(h)=\int c_1(h)^2+c_2(h)$; 
\item[-] 
$D=\int tp(A_3)(g)=\int c_1(g)^3+3c_1(g)c_2(g)+2c_3(g)$. 
\end{itemize}
Computing the degrees in the right hand sides of above equalities 
(i.e. computing Gysin images via $i_*$, $j_*$ and $f_*$),  
we see that  
\begin{eqnarray*}
	m_3&=&4 d - \xi_{001} + 2 \xi_{01} - 3 \xi_1\\
	&=&d(d-1)^3 + (-6  + 9 d  - 4 d^2  + 2 \mu_0)\mu_0 - 4 q \\
	&&\textstyle 
	- \frac{3}{2} s_t  + (4 d-9) t +(\frac{3}{2} d-14) \gamma + 4 \chi_C
\end{eqnarray*}

\end{exam}

\begin{exam} \upshape \label{lem_classD} 
\upshape 
We explain how to express in terms of $d, \cdots, \xi_{001}$ 
numerical characters associated to $D$,  
e.g. rank $\mu_1$ and class $\mu_2$ 
(=rank and class of the image surface $\pi(D) \subset \Proj^3$ under a generic projection). 
Recall that 
the source double point locus $\Gamma \, (\subset M)$ is a surface with ordinary singularities. 
Take a resolution  $\Gamma'$ of $\Gamma$; 
the composed map $\varphi: \Gamma' \to M$ into the ambient $3$-fold  
is a stable map so that 
\begin{equation}\label{vf}
\varphi(A_0^k(\varphi))=A_0^{k+1}(f) \quad (k=1,2,3), \quad 
\varphi(A_1(\varphi))=A_0A_1(f). 
\end{equation}

Let $\mathcal{R}$ be 
the $\Q$-subalgebra of $H^*(M; \Q)$ generated by 
Chern classes $c_i(TM)$,  
the pushforward of all Chern monomials  
$f^*f_*(c_I(TM))$  and $\tilde{a}:=f^*a$. 
Notice that the degree of the image via $f_*$ of an element of $\mathcal{R}$ 
is a polynomial of  $d, \xi_1, \cdots, \xi_{001}$. 
Obviously, $c_i(f)$ and $f^*s_I(f)=f^*f_*(c_I(f))$ belong to $\mathcal{R}$, thus 
universal polynomials $tp$ and $tp^{\SM}$ applied to $f$ are all in $\mathcal{R}$.  
It follows from (\ref{vf}) that 
$$\textstyle 
\frac{1}{k}\varphi_*tp(A_0^k)(\varphi)=tp(A_0^{k+1})(f), \quad 
\varphi_*tp(A_1)(\varphi)=tp(A_0A_1)(f), 
$$
thus those classes are in $\mathcal{R}$. 
Recall that 
\begin{eqnarray*}
\varphi_*tp(A_0)(\varphi)
&=&\textstyle 
\varphi_*(1)=s_0(\varphi), \\
\textstyle 
\varphi_*tp(A_0^2)(\varphi)
&=&\textstyle 
s_0(\varphi)^2-s_1(\varphi), 
\\  
\textstyle 
\varphi_*tp(A_0^3)(\varphi)
&=&\textstyle 
\frac{1}{2}(s_0(\varphi)^3-3s_0(\varphi)s_1(\varphi)+2s_2(\varphi)+2s_{01}(\varphi)),\\
\textstyle
\varphi_*tp(A_1)(\varphi)
&=& \textstyle 
s_{01}(\varphi), 
\end{eqnarray*}
and thus 
\begin{eqnarray*}
s_{0}(\varphi)&=&\varphi_*(1), \\
s_1(\varphi)&=&
c_1(M)\varphi_*(1)-\varphi_*(c_1(\Gamma')), \\
s_2(\varphi)&=&
c_1(M)^2\varphi_*(1)-2c_1(M)\varphi_*(c_1(\Gamma')) +\varphi_*(c_1(\Gamma')^2), \\
s_{01}(\varphi)&=&
c_2(M)\varphi_*(1)-c_1(M)\varphi_*(c_1(\Gamma')) +\varphi_*(c_1(\Gamma')^2)-\varphi_*(c_2(\Gamma')) 
\end{eqnarray*}
are in $\mathcal{R}$.  Hence   
$$\varphi_*(1), \;\; \varphi_*(c_1(T\Gamma')), \;\; 
\varphi_*(c_1(T\Gamma')^2), \;\; 
\varphi_*(c_2(T\Gamma')) \in \mathcal{R}.$$ 
In particular, their degrees  
are written in terms of $d, \cdots, \xi_{001}$. 
To compute numerical characters associated to $D$, 
we take composed maps $\pi\circ f\circ \varphi$ 
with generic linear projections $\pi$ from a point and a line. 
Notice that $\Gamma \to D$ is doubly ramified along smooth $C$. 
In a quite similar way  
as in the computations of the rank and the class 
of a surface in \S 3.2,  
$\mu_1$ and $\mu_2$ of $D$ 
can be interpreted in terms of 
degrees of the critical loci of $\pi\circ f\circ \varphi$, 
its cuspidal locus and $C$. 
Then we see that 
they are expressed by the degree $\mu_0$ of $D$, 
 $\int c_1(T\Gamma')$ and $\int c_2(T\Gamma')$, 
thus by $d, \cdots, \xi_{001}$. 
Apparent characters of $D$ (characters of singular loci of projections of $D$) are also 
computed in the same way. 
\end{exam}

\begin{rem}\upshape \label{tsuboi}
As an analogy to Lemma \ref{enriques}, 
Chern numbers $\xi_3, \xi_{11}, \xi_{001}$ of the normalization $M$ of 
a projective $3$-fold $X$ with ordinary singularities 
have been studied  
 by Tsuboi \cite{Tsuboi05, Tsuboi09} 
in a completely different method using 
the excess intersection formula and Piene's formulae of polar classes. 
Together with seven characters $d, \mu_0, \cdots, \chi_C$, 
his result involves 
the intersection product of the canonical divisor of $M$ with $S=\overline{A_1(f)}\, (\simeq C=f(S))$; 
the number is computed by our method as 
\begin{eqnarray*} 
K_M \cdot S&=&\textstyle -c_1(TM)\cdot tp(A_1)\\
&=&\textstyle 
 -10 \xi_1 + \xi_{11} + 5 \xi_2 - \xi_3\\
&=&\textstyle
\frac{1}{2}s_t - \frac{1}{2}d \gamma - \chi_C 
\end{eqnarray*}
(also note that  $\chi(C, \Ost_C)=\frac{1}{2}\chi_C$ 
for the smooth curve $C$ of crosscap points). 
Then it turns out that 
Tsuboi's formulae in \cite{Tsuboi09} and ours in Table \ref{chern_number} completely 
coincide. 
We also confirm the class $m_3$ of $X$ computed in \cite{Tsuboi09} 
(cf. Example \ref{classX} above). 
\end{rem}

\begin{exam}\upshape \label{tsuboi_ex} (Example 3.3 in \cite{Tsuboi09}) 
Let $\iota: M=\Proj^3 \to \Proj^9$ be the quadratic Veronese embedding:  
$$\iota[z_0; z_1; z_2; z_3] 
=[z_0^2; z_1^2; z_2^2; z_3^2; z_0z_1; z_0z_2; z_0z_3; z_1z_2; z_1z_3; z_2z_3].$$
Let $X$ be the image of $\iota(M)$ via a generic projection to $\Proj^4$, and 
$f: M \to X \subset \Proj^4$ the obtained stable map (Mather \cite{Mather73}). 
Since $c(TM)=(1+a)^4 \in H^*(M)=\Q[a]/\langle a^4\rangle$ 
and the pullback of the hyperplane class is $2a$, 
thus $d = 8$, $\xi_1 = 16$, $\xi_2 = 32$, $\xi_3 = 64$, $\xi_{01} = 12$, $\xi_{11} = 24$, 
$\xi_{001} = 4$. By Table \ref{tp_computation} and Example \ref{classX} above, we have 
$$\mu_0=20, \; t=20, \;\gamma=20, \;q=5, \;s_t=40, \;\chi_C=-20, \; m_3=4.$$
\end{exam}

\subsection*{Acknowledgement} 
This is based on part of the first author's diploma (March 2015, Hokkaido Univerisity). 
The second author thanks Shoji Tsuboi for valuable comments. 
This paper was mostly written during the second author's stay 
in University of Vienna; he thanks for their hospitality. 
This work was partly supported by JSPS KAKENHI Grant Numbers 24340007 and 15K13452.



\begin{thebibliography}{99999}
%
\bibitem{Ando}Y.~Ando, 
On the Higher Thom polynomials of Morin Singularities, 
Publ. RIMS, Kyoto Univ., {\bf 23} (1987), 195--207.
%
\bibitem{Baker} H.~F.~Baker, 
{\it Principles of Geometry, VI. Introduction to the theory of algebraic surfaces and higher loci}, 
Cambridge University Press (1933). 
%
\bibitem{BruceKirk} J.~W. Bruce and N.~Kirk, 
Generic projections of stable mappings, 
{\it Bull. London Math. Soc.} {\bf 32} (2000), 718--728. 
%
\bibitem{Colley} S. ~J. ~Colley, 
Enumerating stationary multiple-points,  
{\it Adv. Math.} {\bf 66} (1987), 149--170. 
%
\bibitem{FR04} L.~Feh\'er and R.~Rim\'anyi, 
Calculation of Thom polynomials 
and other cohomological obstructions for group actions, 
{\it Real and Complex Singularities} (Sao Carlos, 2002), 
Contemp. Math.  {\bf 354}, A.M.S. (2004), 69--93.
%
\bibitem{FR12} L.~Feh\'er and R.~Rim\'anyi, 
Thom series of contact singularities, 
Ann.  Math. {\bf 176} (2012), 1381-1426. 
%
\bibitem{Kaz03} M.~E.~Kazarian,  
Multisingularities, cobordisms and enumerative geometry, 
Russian Math. Survey {\bf 58} (2003), 665--724 (Uspekhi Mat. nauk {\bf 58}, 29--88). 
%
\bibitem{Kaz08} M.~E.~Kazarian, 
Morin maps and their characteristic classes (2008), preprint. 
%
\bibitem{Kleiman76} S.~L.~Kleiman,  
The enumerative theory of singularities, 
Proc. Nordic Summer School/NAVF, Symposium in Math., Oslo, 
Sijthoff  and Noordhoff Inter. Publ. (1976), 297--396. 
%
\bibitem{Kleiman81} S.~L.~Kleiman,  
Multiple-point formulas. I: Iteration, 
Acta Math. {\bf 147} (1981), 13--49.
%
%
\bibitem{Mac} R.~MacPherson, 
 Chern classes for singular algebraic varieties,
{\it Ann. Math.} {\bf 100}  (1974),  421--432.
%
\bibitem{Mather73} J.~Mather, 
Generic projections, 
Ann. Math. {\bf 98} (1973), 226--245. 
%
\bibitem{Mather78} J.~Mather, 
Stable map-germs and algebraic geometry, 
Lecture Notes in Math. {\bf  678}, Springer, (1978), 196--235. 
%
\bibitem{Ohmoto06}T.~Ohmoto, 
Equivariant Chern classes of singular algebraic varieties with group actions,  
{\it Math. Proc. Cambridge Phil. Soc.} {\bf 140} (2006), 115--134.
%
\bibitem{Ohmoto07}T.~Ohmoto, 
Chern classes and Thom polynomials,  
{\it Singularities in Geometry and Topology} (ICTP, Trieste, Italy,  2005), 
World Scientific (2007),  464--482. 
%
\bibitem{Ohmoto16}T.~Ohmoto, 
Singularities of maps and characteristic classes,  
{\it School on Real and Complex Singularities in S\~ao Carlos, 2012}, 
Adv. Studies. Pure Math. {\bf 68} (2016), 191--265. 
ArXiv:1309.0661. 
%
\bibitem{PP}A.~Parusi\'nski and P.~Pragacz, 
Chern-Schwartz-MacPherson classes 
and the Euler characteristic of degeneracy loci and special divisors,
Jour. Amer. Math. Soc. {\bf 8} (1995), 793--817.
%
\bibitem{Piene} R.~Piene, 
Some formulae for a surface in $\Proj^3$, 
{\it Algebraic Geometry}, 
Lecture Notes Math. {\bf 687}, Springer  (1978), 196--235. 
%
%
\bibitem{Roth} L.~Roth, 
Some formulae for primals in four dimensions (I), (II) and (III) 
Proc. London Math. Soc.  (1933), 540--550; 
Proc. London Math. Soc.  (1933), 334--338; 
Proc. Cambridge Phil. Soc.  (1936), 365--369. 
%
\bibitem{Salmon} G.~Salmon, 
{\it A treatise on the analytic geometry of three dimensions},  
4th edition, Dublin (1882). 
%
\bibitem{SO} T.~Sasajima and T.~Ohmoto, 
Thom polynomials in $\A$-classification I: 
counting singular projections of a surface, 
 IMPANGA Lecture Notes `Vector bundles, Schubert varieties, and equivariant cohomology', 
 Birkhauser (2017). arXiv:1606.09147. 
%
\bibitem{SR} J.~G.~Semple and L.~Roth, 
{\it Introduction to Algebraic Geometry},  
Oxford University Press (1949). 
%
\bibitem{Tsuboi05} S.~Tsuboi, 
The Chern numbers of the normalization of an algebraic threefold with ordinary singularities, 
The series ``Seminaires et Congress" 10, Sci. Math. France, (2005), 351--372. 
%
\bibitem{Tsuboi09} S.~Tsuboi, 
Linear projections of rational threefolds, 
Proceedings of the 16th International Conference 
on Finite or Infinite Dimensional Complex Analysis and Applications, 
Daeyang Printing (Gyeongju, Korea), (2009), 
237--247. 
%
%
\end{thebibliography}
\end{document}